\documentclass[12pt, reqno]{amsart}
\usepackage{hyperref}

\usepackage{amssymb,amsmath,graphicx,amsthm}

\def\simge{\underset\sim>}

\def\T{\text}

\def\1#1{\overline{#1}}

\def\2#1{\widetilde{#1}}

\def\3#1{\widehat{#1}}

\def\4#1{\mathbb{#1}}

\def\5#1{\frak{#1}}

\def\6#1{{\mathcal{#1}}}

\def\C{{\4C}}

\def\R{{\4R}}

\def\La{\Lambda}

\def\T{\text}
\newcommand{\Om}{\Omega}
\newcommand{\om}{\omega}

\newcommand{\no}[1]{\|{#1}\|}

\def\R{{\Bbb R}}

\def\C{{\Bbb C}}

\def\di{\partial}
\def\dib{\bar\partial}
\def\Label#1{\label{#1}}
\def\simge{\underset\sim>}

\def\T{\text}

\def\1#1{\overline{#1}}

\def\2#1{\widetilde{#1}}

\def\3#1{\widehat{#1}}

\def\4#1{\mathbb{#1}}

\def\5#1{\frak{#1}}

\def\6#1{{\mathcal{#1}}}

\def\C{{\4C}}

\def\R{{\4R}}

\def\La{\Lambda}

\numberwithin{equation}{section}
\renewcommand\Re{\operatorname{Re}}
\renewcommand\Im{\operatorname{Im}}

\def\T{\text}

\theoremstyle{plain}

\newtheorem{example}{Example}[section]

\newtheorem{theorem}{Theorem}[section]

\newtheorem{lemma}[theorem]{Lemma}

\theoremstyle{definition}

\theoremstyle{remark}

\newtheorem{remark}[theorem]{Remark}

\begin{document}

\title[Supnorm and $f$-H\"{o}lder estimates]{Supnorm and $f$-H\"{o}lder estimates for $\dib$ on convex domains of general type in $\C^2$}
  \author[T.V. Khanh]{Tran Vu Khanh}   
%\address{Tan Tao University, Tan Tao avenue, Tan Duc e-city, Long An province, Vietnam}
%\email{khanh.tran@ttu.edu.vn}
%\address{Dipartimento di Matematica, Universit\`a di Padova, via 
%Trieste 63, 35121 Padova, Italy}
%\email{khanh@math.unipd.it}
\begin{abstract}In this paper, we study supnorm and  modified H\"{o}lder estimates for the integral solution of the $\dib$-equation on a class of convex domains of general type in $\C^2$ that includes many infinite type examples. \\
%AMS Mathematics Subject Classification
%(2000) and key words and phrases.
\\[2mm] {\it AMS Mathematics  Subject Classification (2000)}: 23A26, 32A30, 32T18, 32T25, 32W05. 
\\[1mm] {\it Key words and phrases}: convex domain; domain of finite/infinite type; $f$-H\"older, supnorm, $L^2$-superlogarithmic estimate.
\end{abstract}

\maketitle\tableofcontents 
\section{Introduction}
   % Let $\Om$ be a bounded pseudoconvex domain in $\C^2$ with smooth boundary and $0\in b\Om$ .
%;  and there is a constant $\beta>1$ such that 
%$$1\le \frac{tF'(t)}{F(t)}\le \frac{1}{t\ln^\beta t}\quad t\in (0,\delta)$$
 %===============

Let $\Om$ be a smooth, bounded domain in $\C^2$ with $0$ in the boundary $b\Om$. Assume that $\Om$ is strictly convex except possibly on a neighborhood $U$ of $0$;  and in $U$,  $\Om$ has the form 
\begin{eqnarray}\Label{Om1}
\Om\cap U=\{\rho(z)=F(|z_1|^2)+r(z)<0\}
\end{eqnarray}
or 
\begin{eqnarray}\Label{Om2}
\Om\cap U=\{\rho(z)=F(|\Re{z_1}|^2)+r(z)<0\}.
\end{eqnarray}
where $F$ is a strictly increasing, convex function such that $F(0)=0$, $F(t)/t$ is increasing, and  $r$ is convex with $\dfrac{\di r}{\di z_2}\not=0$. We remark that $\Om$ may be of finite type or infinite type since we may choose, for example, $F(t)=t^m$ or $F(t) = \exp(-1/t^\alpha)$.
%We remark that $\Om$ has the general type which is smaller than or equal to $F$ (including finite or infinite type, for example: $F(t)=t^m$ or $F(t)=\exp(-1/t^\alpha)$). 
The primary goals of this paper are to investigate the supnorm estimate and develop appropriate H\"older estimate for the integral solution of the $\dib$-equation given by Henkin kernel on a domain $\Om$ satisfying \eqref{Om1} or \eqref{Om2}. \\

Given a  bounded, $\dib$-closed $(0,1)$ form $\phi$, the supnorm and the H\"older estimates for the solution of the Cauchy-Riemann equation 
$$\dib u=\phi$$  on domain $\Om$ is a fundamental question in several complex variables. A positive answer is well-known when $\Om$ is a
\begin{enumerate}
  \item[$\bullet$]  strongly pseudoconvex domain in $\C^n$ (see \cite{He70}, \cite{Ke71}, \cite{Ra86} ...),
  \item[$\bullet$]  convex domain of finite type in $\C^n$ (see \cite{DF06}, \cite{DFF99}, \cite{H02}...),
  \item[$\bullet$]  real or complex ellipsoid of finite type in $\C^n$ (see \cite{BC84}, \cite{F96}, \cite{DFW86},... ),
  \item[$\bullet$]  or a pseudoconvex domain of finite type in $\C^2$ (see \cite{FK88},\cite{Ra90},\cite{CNS92}...).
  \end{enumerate} 
However, when $\Omega$ is of infinite type, the only result is by  J.~E.~Fornaess, L.~Lee, and Y.~Zhang \cite{FLZ11} who prove supnorm estimates in the case $F(t) = \exp(-1/t^\alpha)$ with $\alpha <\dfrac{1}{2}$ and $r(z)=\Re{z_2}$  for both \eqref{Om1} and \eqref{Om2}. Denote by $A\lesssim B$ for inequality $A\le cB$ with some positive constant $c$, for simplification. We denote by $L^\infty(\Om)$ the space of  the essentially bounded functions on $\Om$ and by $\no{u}_{\infty}$ the essential supremun of $u\in L^\infty(\Om)$ in $\Om$.
\begin{theorem}[Fornaess-Lee-Zhang]
Let $\Om$ be a smooth, bounded domain in $\C^2$ with $0$ in the boundary $b\Om$. Assume that $\Om$ is strictly convex except possibly on a neighborhood $U$ of $0$;  and in $U$,  $\Om$ has the form 
\begin{eqnarray}
\Om\cap U=\{\rho(z)=\Re{z_2}+\exp(-1/|z_1|^\alpha)<0\}
\end{eqnarray}
or 
\begin{eqnarray}
\Om\cap U=\{\rho(z)=\Re{z_2}+\exp(-1/|\Re z_1|^\alpha)<0\}.
\end{eqnarray}
with $\alpha<1$. Then there is a solution to the $\dib$-equation $\dib u=\phi$ for any $\phi\in C^1_{(0,1)}(\bar\Om)$ and $\dib \phi=0$, that satisfies $\no{u}_\infty\lesssim \no{\phi}_\infty$.

\end{theorem}
% The first goal of this paper is to give the proof of supnorm estimates which is slightly more general than the result in \cite{FLZ11}. 
The first goal of the paper is to prove supnorm estimates on domains satisfying \eqref{Om1} or \eqref{Om2} which both generalize the class of domains considered in \cite{FLZ11}. 
%and decrease the allowable smoothness of $\phi$ from $C^1_{(0,1)}(\bar\Omega)$ to $(L_\infty)_{(0,1)}(\bar\Omega)$.

\begin{theorem}\Label{t1} (i) Let $\Om$ and $F$ be as in \eqref{Om1}. Assume that $\displaystyle\int_0^\delta |\ln F(t^2)|dt<\infty$ for some $\delta>0$. Then for any bounded $\dib$-closed $(0,1)$ form $\phi$ on $\bar\Om$, there is a $u$ such that $\dib u=\phi$ on $\Om$ and $\no{u}_\infty\lesssim \no{\phi}_\infty$.  \\

(ii) Let $\Om$ and $F$ be as in \eqref{Om2}. Assume that $\displaystyle\int_0^\delta |(\ln t)(\ln F(t^2))|dt<\infty$ for some $\delta>0$. Then for any bounded $\dib$-closed $(0,1)$ form $\phi$ on $\bar\Om$, there is a $u$ such that $\dib u=\phi$ on $\Om$ and $\no{u}_\infty\lesssim \no{\phi}_\infty$.  
\end{theorem}

When $\Om$ is finite type (e.g., $F(t)=t^m$),  we known that fractional H\"older estimates hold  for both case \eqref{Om1} and \eqref{Om2}. However, when $\Om$ is infinite type (e.g., $F(t)=\exp(-1/t^\alpha)$), McNeal \cite{Mc91} proves the fractional H\"{o}lder estimates do not hold. \\

In this paper, we find a suitable H\"older estimate for infinite type. Let $f$ be an increasing function on $(a,+\infty)$ with $a$ big enough such that $\underset{t\to+\infty}{\lim} f(t)=+\infty$. For $\Om\subset \C^n$, we define $f$-H\"older space on $\Om$ by
$$\Lambda^{f}(\Om)=\{u : \no{u}_{\infty}+\sup_{z,z+h\in \Om}f(|h|^{-1}) \cdot |u(z+h)-u(z)|<\infty \}$$ 
and set 
$$\no{u}_{f}= \no{u}_{\infty}+\sup_{z,z+h\in \Om}f(|h|^{-1})\cdot |u(z+h)-u(z)|. $$
Note that the $f$-H\"older spaces include the standard H\"older spaces $\Lambda_\alpha(U)$ by taking $f(t) = t^{\alpha}$ (so $f(|h|^{-1}) = |h|^{-\alpha}$) with $0<\alpha<1$. In this way, $f$-H\"older spaces generalize the notion of the H\"older spaces.\\

Since $F$ is strictly increasing $F$ is invertible with inverse $F^*$. Our main result is
\begin{theorem}\Label{t2}(i) Let $\Om$ and $F$ be defined by \eqref{Om1}.  Assume that $\displaystyle\int_0^\delta |\ln F(t^2)|dt<\infty$ for some $\delta>0$.  Then for every bounded $\dib$-closed $(0,1)$ form $\phi$ on $\bar\Om$, there exists a function $u$ in $\Lambda^{f}(\Om)$ such that $\dib u=\phi$ and 
\begin{eqnarray}
\no{u}_{f} \lesssim\no{\phi}_{\infty}
\end{eqnarray}
where $$f(d^{-1})=\Big(\displaystyle\int_0^d\frac{\sqrt{F^*(t)}}{t}dt\Big)^{-1}.$$ 

(ii) Let $\Om$ and $F$ be defined by \eqref{Om2}.  Assume that $\displaystyle\int_0^\delta |(\ln t) (\ln F(t^2))|dt<\infty$ for some $\delta>0$.  Then for every bounded $\dib$-closed $(0,1)$ form $\phi$ on $\bar\Om$, there exists a function $u$ in $\Lambda^{f}(\Om)$ such that $\dib u=\phi$ and 
\begin{eqnarray}
\no{u}_{f} \lesssim \no{\phi}_{\infty}
\end{eqnarray}
where $$f(d^{-1})=\Big(\displaystyle\int_0^d\frac{\sqrt{F^*(t)}\big|\ln\sqrt{F^*(t)}\big|}{t}dt\Big)^{-1}.$$ \\
\end{theorem}
The following examples are explicit function $f$ in the choice of $F$.
\begin{example}{\rm  Let $\Om=\{(z_1,z_2)\in \C^2: F(|z_1|^2)+|z_2-1|^2<1\}$. Then  supnorm and $f$-H\"older estimates hold for the integral solution of $\dib$-equation in the following examples:
\begin{enumerate}
  \item If $F(t^2)=t^{2m}$, then $f(d^{-1})=d^{-1/2m}$.
  \item If $F(t^2)=2\exp\left(-\frac{1}{t^\alpha}\right)$ with $0<\alpha<1$, then $f(d^{-1})=(-\ln d)^{\frac{1}{\alpha}-1}$.
  \item  If $F(t^2)=2\exp\left(-\frac{1}{t|\ln t|^\alpha}\right)$ with $\alpha>1$, then $f(d^{-1})=\big(\ln(-\ln d)\big)^{\alpha-1}$
\end{enumerate}}
\end{example}

\begin{example}{\rm  Let $\Om=\{(z_1,z_2)\in \C^2: F(|\Re z_1|^2)+G(|\Im z_1|^2)+|z_2-1|^2<1\}$ where $G(t)\equiv0$ in a neighborhood of $0$ and there is a positive constant $c$ such that $t\ge c$ if $G(t)\ge 1$. Then supnorm and $f$-H\"older estimates hold for the integral solution of the $\dib$-equation in the following examples:
\begin{enumerate}
  \item If $F(t^2)=t^{2m}$, then $f(d^{-1})=d^{-1/2m}|\ln d|^{-1}$.
  \item If $F(t^2)=\exp\left(-\frac{1}{t^\alpha}\right)$ with $0<\alpha<1$, then $f(d^{-1})=(-\ln d)^{\frac{1}{\alpha}-1}\big(\ln(-\ln d)\big)^{-1}$.
  \item  If $F(t^2)=\exp\left(-\frac{1}{t|\ln t|^\alpha}\right)$ with $\alpha>2$, then $f(d^{-1})=\big(\ln(-\ln d)\big)^{\alpha-2}$.
\end{enumerate}}
\end{example}
\begin{remark}\Label{r1} We remark that superlogarithmic estimates defined by Kohn in \cite{Ko02} for the $\dib$-Laplacian $\Box$ or Kohn-Laplacian $\Box_b$, which imply local hypoellipcity of $\Box$ and $\Box_b$ respectively,  hold in both domains defined by \eqref{Om1} and $\eqref{Om2}$ under hypothesis in Theorem \ref{t2} (see Appendix).   
\end{remark}

I am grateful to Andrew Raich for helpful comments on the original draft of this paper.
\section{Preliminaries}
In this section we briefly recall  the contruction of an integral kernel to  solve the $\dib$-equation on convex domains in $\C^2$. Details can be found in \cite{He70}, \cite{Ra86}.  \\
 
Let $\rho$ be the defining function of $\Om$. We can assume that there is a $\delta>0$ such that $b\Om\setminus B(0,\delta)$ is strictly convex and 
\begin{eqnarray}\Label{Om}
\Om\cap B(0,\delta)=\{\rho(z)=P(z_1)+r(z)<0\}
\end{eqnarray}
where $P(z_1)=F(|z_1|^2)$ or $F(|\Re z_1|^2)$ and $r(z)$ is convex with $\dfrac{\di r}{\di z_2}\not=0$ on $\Om\cap B(0,\delta)$.

Define 
$$\Phi(z,\zeta)=2\Big(\frac{\di \rho(\zeta)}{\di \zeta_1}(\zeta_1-z_1)+\frac{\di \rho(\zeta)}{\di \zeta_2}(\zeta_2-z_2)\Big).$$
The following result is a well-known consequence of Taylor's theorem and convexity
\begin{lemma}\Label{l2.1}
There are suitably small $\epsilon$ and $c$ such that 
\begin{eqnarray}\Label{ReQ}
\Re \Phi(z,\zeta)\ge -\rho(z)+\begin{cases} c|z-\zeta|^2  & \zeta\in b\Om\setminus B(0,\delta)\\
                                                                         P(z_1)-P(\zeta_1)-2\Re\dfrac{\di P}{\di \zeta_1}(\zeta_1)(z_1-\zeta_1) & \zeta\in b\Om\cap B(0,\delta) \end{cases}
\end{eqnarray}
for all $z\in \bar\Om$ with $|z-\zeta|\le \epsilon$. 
\end{lemma}  
Choose $\chi\in C^\infty(\C^2\times \C^2)$ such that $0\le \chi\le 1$, $\chi(z,\zeta)=1$ for $|z-\zeta|\le \dfrac{1}{2}\epsilon$ and $\chi(z,\zeta)=0$ for $|z-\zeta|\ge \epsilon$. For $j=1,2$, define 
$$\Phi_j(z, \zeta)=\chi\frac{\di \rho}{\di \zeta_j}(\zeta)+(1-\chi)(\bar\zeta_j-\bar z_j).$$
Then
$$\Phi^\#(z,\zeta)=\Phi_1(z,\zeta)(\zeta_1-z_1)+\Phi_2(z,\zeta)(\zeta_2-z_2)$$
has the following properties for any $\zeta\in b\Om$:
\begin{enumerate}
  \item[(i)] \begin{eqnarray}\Re\Phi^\#(z,\zeta)\ge  -\rho(z)+\begin{cases} c|z-\zeta|^2  & \zeta\in b\Om\setminus B(0,\delta),\\
                                                                         P(z_1)-P(\zeta_1)-2\Re\dfrac{\di P}{\di \zeta_1}(\zeta_1)(z_1-\zeta_1) & \zeta\in b\Om\cap B(0,\delta); \end{cases}
\end{eqnarray}
for $|z-\zeta|\le \dfrac{1}{2}\epsilon$ and $z\in\bar\Om$.
 % \item[(ii)] there is $\delta>0$ such that $|\Phi(z,\zeta)|\ge \delta$ 
  %for all $z$ with $r(z)\le -\delta$ and $|z-\zeta|\ge \epsilon$; 
\item[(ii)] $\Phi^\#(\cdot, \zeta)$ and $\Phi_j(\cdot, \zeta)$, $j=1,2$, are holomorphic on $\{z: |z-\zeta|\le \dfrac{1}{2}\epsilon\}$.
  \end{enumerate}

We now are ready for the integral solution of the $\dib$-equation. 
Let $\phi=\phi_1d\bar z_1+\phi_2d\bar z_2$ be a bounded $\dib$-closed $(0,1)$-form on $\bar\Om$. 
The Hekin integral solution $u$ of the $\dib$-equation $\dib u=\phi$ given by
 $$u=T\phi(z)=H\phi(z)+K\phi(z).$$
where
\begin{eqnarray}\begin{split}
H\phi(z)=& -\frac{1}{2\pi i}\int_{\zeta\in b\Om}\frac{\Phi_1(z, \zeta)(\bar\zeta_2-\bar z_2)-\Phi_2(z, \zeta)(\bar\zeta_1-\bar z_1)}{\Phi^\#(z,\zeta)|\zeta-z|^2}\phi\wedge\om(\zeta);\\
K\phi(z)=&-\frac{1}{2\pi i}\int_\Om \frac{\phi_1(\bar\zeta_2-\bar z_2)-\phi_2(\bar\zeta_1-\bar z_1)}{|\zeta-z|^4}\om(\bar\zeta)\wedge\om(\zeta)
\end{split}
\end{eqnarray}
where $\om(\zeta)=d\zeta_1 \wedge d\zeta_2$.\\

It is well known that 
$$\no{K\phi}_\infty\lesssim \no{\phi}_\infty \quad\T{and}\quad \no{K\phi}_f\lesssim  \no{\phi}_\infty,$$
for any $0<f(d^{-1})<d^{-1}$ with any $d>0$ small enough (see Lemma 1.15, page 157 in \cite{Ra86}). Moreover, we have 
$$H\phi(z)=-\frac{1}{2\pi i}\int_{\zeta\in b\Om}\cdots=-\frac{1}{2\pi i}\Big(\int_{\zeta\in b\Om, |z-\zeta|\le \epsilon}\cdots+\int_{\zeta\in b\Om, |z-\zeta|\ge\epsilon}\cdots\Big)=-\frac{1}{2\pi i}\int_{\zeta\in b\Om, |z-\zeta|\le \epsilon}\cdots$$
Since $\Phi_1(z, \zeta)(\bar\zeta_2-\bar z_2)-\Phi_2(z, \zeta)(\bar\zeta_1-\bar z_1)\equiv0$ if $|z-\zeta|\le \epsilon$.\\
Therefore, it is sufficient to estimate 
$$H\phi(z)=-\frac{1}{2\pi i}\int_{\zeta\in b\Om, |z-\zeta|\le \epsilon}\frac{\Phi_1(z, \zeta)(\bar\zeta_2-\bar z_2)-\Phi_2(z, \zeta)(\bar\zeta_1-\bar z_1)}{\Phi(z,\zeta)|\zeta-z|^2}\phi\wedge\om(\zeta).$$

We will use the general Hardy-Littewood lemma (see Section 5 below) to obtain the $f$-H\"older estimates. To do that we need to control the gradient of $T\phi(z)$. 
We have 
$$|\nabla H\phi(z)|\lesssim \no{\phi}_\infty\int_{\zeta\in b\Om, |z-\zeta|\le \epsilon}\Big(\dfrac{1}{|\Phi|\cdot|\zeta-z|^2}+\dfrac{1}{|\Phi|^2\cdot|\zeta-z|}\Big)dS$$
where $dS$ is surface area measure on $b\Om$. 
%For $a>0$, to be determined later, let $B$ denote the ball with center $(0,1)$ and radius $a$. 
We now use Lemma \ref{l2.1} to obtain
$$\int_{\zeta\in b\Om\setminus B(0,\delta), |z-\zeta|\le \epsilon}\Big(\dfrac{1}{|\Phi|\cdot|\zeta-z|^2}+\dfrac{1}{|\Phi|^2\cdot|\zeta-z|}\Big)dS\lesssim \delta^{-1/2}(z).$$
Hence, it remains to estimate 
$$L(z)=\int_{\zeta\in b\Om\cap B(0,\delta), |z-\zeta|\le \epsilon}\Big(\dfrac{1}{|\Phi|\cdot|\zeta-z|^2}+\dfrac{1}{|\Phi|^2\cdot|\zeta-z|}\Big)dS.$$
Set $t=\Im\Phi(z, \zeta)$. It is easy to check that $\dfrac{\di t}{\di \zeta_2}\not=0$. So we change coordinate and  obtain
\begin{eqnarray}\begin{split}
L(z)\lesssim &\int_{|t|\le \delta, |\zeta_1|<\delta,|z_1-\zeta_1|\le \epsilon }\frac{dtd(\Re\zeta_1)d(\Im\zeta_1)}{(|t|+|\Re\Phi|)(|\rho(z)|^2+|\zeta_1-z_1|^2)}\\
&+\int_{|t|\le \delta, |\zeta_1|<\delta,|z_1-\zeta_1|\le \epsilon } \frac{dtd(\Re\zeta_1)d(\Im\zeta_1)}{(t^2+|\Re\Phi|^2)(|\rho(z)|+|\zeta_1-z_1|)}\\
 \lesssim& |\ln(|\Re\Phi|)\cdot \ln\rho(z)|+ \int_{|\zeta_1|<\delta, |z_1-\zeta_1|\le \epsilon}\frac{d(\Re\zeta_1)d(\Im\zeta_1)}{|\Re\Phi| \cdot|\zeta_1-z_1|}\\
 \lesssim& |\ln\rho(z)|^2+ \int_{|\zeta_1|<\delta, |z_1-\zeta_1|\le \epsilon}\frac{d(\Re\zeta_1)d(\Im\zeta_1)}{|\Re\Phi| \cdot |\zeta_1-z_1|}.\\
\end{split}
\end{eqnarray}
Here the last inequality follows by $|\Re\Phi|\ge |\rho(z)|$ for all $\zeta\in b\Om\cap B(0,\delta)$ and $|z-\zeta|\le\epsilon $ which is itself a consequence of Lemma \ref{l2.1} and the convexity of $P$ (see \eqref{3.3} and \eqref{3.1a} below). \\

We have therefore show 
\begin{eqnarray}\Label{H}
|\nabla H\phi(z)|\lesssim \Big(\rho(z)^{-1/2}+\int_{|\zeta_1|<\delta, |\zeta_1-z_1|<\epsilon}\frac{d(\Re\zeta_1) d(\Im\zeta_1)}{|\Re\Phi|\cdot|\zeta_1-z_1|}\Big)\no{\phi}_\infty.
\end{eqnarray}
A similar argument also shows
\begin{eqnarray}\Label{S}
|H\phi(z)|\lesssim \Big(1+\int_{|\zeta_1|<\delta, |\zeta_1-z_1|<\epsilon}\frac{\big|\ln|\T{Re }\Phi |\big|d(\Re\zeta_1) d(\Im\zeta_1)}{|\zeta_1-z_1|}\Big)\no{\phi}_\infty.
\end{eqnarray}

\section{Estimates on $\Om\cap U=\{\rho(z)=F(|z_1|^2)+r(z)<0\}$}
In this section, we give the proof of Theorem \ref{t1}.(i) and Theorem \ref{t2}.(ii). It is sufficient to estimate the integrals in \eqref{H} and \eqref{S} when $z\in B(0,\delta)$ so the defining function $\rho$ is of the form $\rho(z)=F(|z_1|^2)+r(z)$ in $B(0,\delta)$.

\begin{lemma}\Label{l1} Let $F$ be a convex, $C^2$-smooth function on $[0,\delta]$. Then we have 
$$F(p)-F(q)-F'(q)(p-q)\ge 0$$
for any $p,q\in [0,\delta]$.  
\end{lemma}
The proof is simple and is omitted here.

\begin{lemma}\Label{l2}For $\delta>0$ small enough, let $F$ be an invertible on $[0,\delta]$ such that $\dfrac{F(t)}{t}$ is increasing  $[0,\delta]$. Then 
$$\int^\delta_0 \frac{dr}{\varrho+F(r^2)}\lesssim \frac{\sqrt{F^*(\varrho)}}{\varrho}$$
for any sufficiently small $\varrho>0$.  
\end{lemma}
{\it Proof.} We split our integration to be two terms
$$\int^\delta_0 \frac{dr}{\varrho+F(r^2)}=\int^{\sqrt{F^*(\varrho)}}_0\cdots+\int_{\sqrt{F^*(\varrho)}}^\delta\cdots$$
For the first term, it is easy to see that
$$\int^{\sqrt{F^*(\varrho)}}_0 \frac{dr}{\varrho+F(r^2)}\le \frac{\sqrt{F^*(\varrho)}}{\varrho}.$$
Since $\dfrac{F(t)}{t}$ is increasing, we have 
$$\frac{F(r^2)}{r^2}\ge \frac{F(F^*(\varrho))}{F^*(\varrho)}=\frac{\varrho}{F^*(\varrho)}, \quad \T{or}\quad  \frac{F(r^2)}{\varrho}\ge \frac{r^2}{F^*(\varrho)},$$
for any $r\ge \sqrt{F^*(\varrho)}$.  Apply this inequality to the second term, we obtain
\begin{eqnarray}\Label{3.2a}\begin{split}
\int_{\sqrt{F^*(\varrho)}}^\delta  \frac{dr}{\varrho+F(r^2)}\le &\frac{1}{\varrho} \int_{\sqrt{F^*(\varrho)}}^\delta \frac{dr}{1+r^2/ F^*(\varrho)}\\
\le & \frac{\sqrt{F^*(\varrho)}}{\varrho}  \int_{1}^\infty \frac{dy}{1+y^2}=\frac{\pi}{4}\frac{\sqrt{F^*(\varrho)}}{\varrho}
\end{split}
\end{eqnarray} 
This is complete the proof of Lemma \ref{l2}.

$\hfill\Box$

{\it Proof of Theorem \ref{t1}.(i).} We omit the proof of Theorem \ref{t1}.(i) since it follows in exactly method of the proof of Theorem \ref{t2}.(i) with simpler calculation.\\    

$\hfill\Box$

{\it Proof of Theorem \ref{t2}.(i). }
We apply the identity $2\Re{a\bar b}=|a+b|^2-|a|^2-|b|^2$ in \eqref{ReQ} to obtain
\begin{eqnarray}\Label{3.3}\begin{split}
\Re\Phi(z,\zeta)\ge& -\rho(z)+F(|z_1|^2)-F(|\zeta_1|^2)+2F'(|\zeta_1|^2)\Re\big(\bar\zeta_1(z_1-\zeta_1)\big).\\
\ge&-\rho(z)+\Big(F'(|\zeta_1|^2)|z_1-\zeta_1|^2+F(|z_1|^2)-F(|\zeta_1|^2)-F'(|\zeta_1|^2)(|z_1|^2-|\zeta_1|^2)\Big)\\
\ge&-\rho(z)+\Big(F'(|\zeta_1|^2)|z_1-\zeta_1|^2)\Big)\\
\end{split}
\end{eqnarray}
where the last inequality follows by Lemma \ref{l1}.\\

Let $M(z)$ be the integral term in \eqref{H}. We will show that 
\begin{eqnarray}\Label{3.4}\begin{split}
M(z)\lesssim \frac{\sqrt{F^*(|\rho(z)|)}}{|\rho(z)|}
\end{split}
\end{eqnarray} 
for $z\in\Om$. For convenient, set $\varrho=|\rho(z)|>0$ when $z\in\Om$. 
From \eqref{3.3}, we have
$$M(z)\le \int_{|\zeta_1|<\delta, |\zeta_1-z_1|<\epsilon}\frac{d(\Re\zeta_1) d(\Im\zeta_1)}{(\varrho+F'(|\zeta_1|^2)|z_1-\zeta_1|^2) |\zeta_1-z_1|}.$$ 
There are now two cases.\\
{\it \underline{Case 1:}} $|z_1-\zeta_1|\ge|\zeta_1|$. In this case,  
$$(\varrho+F'(|\zeta_1|^2)|z_1-\zeta_1|^2)|z_1-\zeta_1|\ge(\varrho+F'(|\zeta_1|^2)|\zeta_1|^2)|\zeta_1|\ge (\varrho+F(|\zeta_1|^2)|\zeta_1|. $$
Here the last inequality follows from the inequality $tF'(t)\ge F(t)$ which is itself a consequence of the fact that $\dfrac{F(t)}{t}$ is increasing.
Therefore, using polar coordinates and Lemma \ref{l2},
\begin{eqnarray}\begin{split}
M(z) \le& \int_{|\zeta_1|<\delta}\frac{d(\Re\zeta_1) d(\Im\zeta_1)}{(\varrho+F(|\zeta_1|^2)) |\zeta_1|}\\
\lesssim&\int_0^\delta \frac{dr }{\varrho+F(r^2)}\\
\lesssim&\frac{\sqrt{F^*(\varrho)}}{\varrho}.
\end{split}
\end{eqnarray}
{\it \underline{Case 2:}} If $|\zeta_1|\ge |z_1-\zeta_1|$, then the fact that $F'$ is increasing ($F$ is convex) implies
 $$F'(|\zeta_1|^2)|z_1-\zeta_1|^2\ge F'(|z_1-\zeta_1|^2)|z_1-\zeta_1|^2\simge F(|z_1-\zeta_1|^2). $$
Similarly, we obtain 
\begin{eqnarray}\begin{split}
M(z)\le& \int_{|\zeta_1-z_1|<\epsilon}\frac{d(\Re\zeta_1) d(\Im\zeta_1)}{(\varrho+F(|\zeta_1-z_1|^2)) |\zeta_1-z_1|}\\
\lesssim&\int_0^\epsilon \frac{dr }{\varrho+F(r^2)}\\
\lesssim & \frac{\sqrt{F^*(\varrho)}}{\varrho}
\end{split}
\end{eqnarray}
The proof of \eqref{3.4} is complete. Combining \eqref{H} and \eqref{3.4}, we obtain 
\begin{eqnarray}
|\nabla T(\phi)|\lesssim \frac{\sqrt{F^*(|\rho(z)|)}}{|\rho(z)|}\no{\phi}_\infty\lesssim \frac{\sqrt{F^*(\delta_{b\Om}(z))}}{\delta_{b\Om}(z)}\no{\phi}_\infty.
\end{eqnarray}
since the distance $\delta_{b\Om}(z)$ is comparable to $|\rho(z)|$. \\

Finally, to apply the general Hardy-Littlewood Lemma (see Section 5), we need to check that $G(t):=\sqrt{F^*(t)}$ satisfies the hypothesis of Theorem \ref{HL}. It is easy to see that $\sqrt{F^*(t)}$ is increasing and $\dfrac{\sqrt{F^*(t)}}{t}$ is decreasing. 
For $\delta$ small enough, $|\ln(F(t^2))|$ is decreasing when $0\leq t \leq \delta$ so we can estimate
\[
\big|\ln F(\eta^2)\big| \eta \leq \int_0^\eta \big|\ln F(t^2)\big|dt \le \int_0^\delta \big|\ln F(t^2)\big|dt <\infty
\]
for any $0\le \eta\le \delta$.
The integral is finite by the hypothesis. Consequently, $\sqrt{F^*(t)}|\ln t| <\infty$ for any $0\leq t \leq \sqrt{F^*(\delta)}$ and $\displaystyle\lim_{t\to0}t|\ln F(t^2)|=0$. This implies
\begin{eqnarray}\Label{3.8}
\begin{split}
\int^d_0\dfrac{\sqrt{F^*(t)}}{t}dt\overset{y:=\sqrt{F^*(t)}}{=}&\int^{\sqrt{F^*(d)}}_0 y\big(\ln F(y^2)\big)' dy\\
=&\sqrt{F^*(d)}\ln d-\int^{\sqrt{F^*(d)}}_0 \big(\ln F(y^2)\big) dy<\infty.
\end{split}
\end{eqnarray} 
for $d$ sufficiently small. Here, the integral in \eqref{3.8} is finite by the hypothesis. Thus the proof of Theorem \ref{t2}.(i) is complete.

$\hfill\Box$
\section{Estimates on $\Om\cap U=\{\rho(z)=F(|\T{\rm Re} z_1|^2)+r(z)<0\}$}
In this section, we give the proof of Theorem \ref{t1}.(ii) and Theorem \ref{t2}.(ii). Is is sufficient to estimate the integrations in \eqref{H} and \eqref{S} when the defining function of $\Om$ in a neigborhood of $0$ has the form $\rho=F(|\Re z_1|^2)+r(z)$. \\

We set $x_1=\T{Re }z_1$, $y_1=\Im z_1$,  $\xi_1=\T{Re }\zeta_1$ and $\eta_1=\Im \zeta_1$. From $\eqref{ReQ}$, we have 
\begin{eqnarray}\Label{3.1a}
\begin{split}
\T{Re }\Phi(z,\zeta)\ge& -\rho(z)+\Big(F(x_1^2)-F(\xi_1^2)-2F'(\xi_1^2)\xi_1(x_1-\xi_1)\Big)\\
\ge& -\rho(z)+F'(\xi_1^2)(x_1-\xi_1)^2+\Big(F(x_1^2)-F(\xi_1^2)-F'(\xi_1^2)(x_1^2-\xi_1^2)\Big)\\
\ge& -\rho(z)+F'(\xi_1^2)(x_1-\xi_1)^2
\end{split}
\end{eqnarray}
where the last inequality follows by Lemma \ref{l1}.\\

{\it Proof of Theorem \ref{t1}.(ii).} We only need to show that the integral term in \eqref{S} is bounded. 
By the estimates of $\T{Re }\Phi(z,\zeta)$ as above, we get 
\begin{eqnarray}\Label{4.2a}
\begin{split}
&\int_{|\zeta_1|<\delta, |\zeta_1-z_1|<\epsilon}\frac{\big|\ln|\T{Re }\Phi |\big|d(\Re\zeta_1) d(\Im\zeta_1)}{|\zeta_1-z_1|}\\
\lesssim&\int_{|\zeta_1|<\delta, |\zeta_1-z_1|<\epsilon}\frac{\big|\ln(F'(\xi_1^2)(x_1-\xi_1)^2)|\big|d\xi_1 d\eta_1}{|x_1-\xi_1|+|y_1-\eta_1|}\\
\lesssim&\int_{|x_1|<\delta, |x_1-\xi_1|<\epsilon}\left|\ln|x_1-\xi_1|\cdot\ln(F'(\xi_1^2)(x_1-\xi_1)^2)|\right|dx_1\\
\lesssim&\int_{|x_1|<\delta, |x_1-\xi_1|<\epsilon; |\xi_1|\le |x_1-\xi_1|}\cdots+\int_{|x_1|<\delta, |x_1-\xi_1|<\epsilon; |\xi_1|\ge |x_1-\xi_1|}\cdots\\
\lesssim&\int_{|x_1|<\delta}\left|\ln|\xi_1|\cdot\ln(F'(\xi_1^2)\xi_1^2)|\right|dx_1\\
~~~~~~&+\int_{|x_1-\xi_1|<\epsilon}\left|\ln|x_1-\xi_1|\cdot\ln(F'((x_1-\xi_1)^2)(x_1-\xi_1)^2)|\right|dx_1\\
\lesssim&\int_{|t|<\max\{\delta,\epsilon\}}\left|\ln|t|\cdot\ln(F(t^2))|\right|dt<\infty
\end{split}
\end{eqnarray}
Here, the first inequality in the last line of \eqref{4.2a} follows from the inequality $t^2F(t^2)\ge F(t^2)$ and the last one of this line follows by hypothesis of theorem. This completes the proof of Theorem \ref{t1}.(ii).

$\hfill\Box$

{\it Proof of Theorem \ref{t2}.(ii).} We only need to estimate of the integral term in \eqref{H}. 
By the estimates of $\T{Re }\Phi(z,\zeta)$ above, we observe
\begin{eqnarray}\Label{4.3a}
\begin{split}
\int_{|\zeta_1|<\delta, |\zeta_1-z_1|<\epsilon}&\frac{d\xi_1 d\eta_1}{(\varrho+F'(\xi_1^2)(x_1-\xi_1)^2)(|x_1-\xi_1|+|y_1-\eta_1|)}\\
\lesssim&\int_{|\xi_1|<\delta, |x_1-\xi_1|<\epsilon}\frac{\big|\ln|x_1-\xi_1|\big|d\xi_1}{\varrho+F'(\xi_1^2)(x_1-\xi_1)^2}\\
\lesssim&\int_{|t|<\max\{\delta, \epsilon\}}\frac{|\ln t|dt}{\varrho+F(t^2)}\\
\end{split}
\end{eqnarray}
Here, the last inequality in the last line of \eqref{4.3a} follows by the comparison of $|\xi_1|$ and $|x_1-\xi_1|$; and the property $t^2F'(t^2)\ge F(t^2)$ as in Theorem \ref{t1}.(ii). To estimate the integral term in the last line of \eqref{4.3a} we need following lemma.
\begin{lemma}\Label{l3}For $\delta>0$ small enough, let $F$ be an invertible on $[0,\delta]$ such that $\dfrac{F(t)}{t}$ is increasing  $[0,\delta]$. Then 
$$\int^\delta_0 \frac{|\ln t|dt}{\varrho+F(t^2)}\lesssim \frac{\sqrt{F^*(\varrho)}|\ln \sqrt{F^*(\varrho)}|}{\varrho}$$
for any $\varrho>0$ sufficiently small.  
\end{lemma}
{\it Proof of Lemma \ref{l3}.}
{\it Proof.} We split our integration into two terms
\begin{eqnarray}\Label{4.5a}
\int^\delta_0 \frac{|\ln t| dt}{\varrho+F(t^2)}=\int^{\sqrt{F^*(\varrho)}}_0\cdots+\int_{\sqrt{F^*(\varrho)}}^\delta\cdots
\end{eqnarray}
For the first term, we have
$$\int^{\sqrt{F^*(\varrho)}}_0 \cdots \le \frac{1}{\varrho}\int_0^{\sqrt{F^*(\varrho)}}|\ln t|dt\lesssim \frac{\sqrt{F^*(\varrho)}|\ln \sqrt{F^*(\varrho)}|}{\varrho}.$$
For the second term 
$$\int_{\sqrt{F^*(\varrho)}}^\delta\cdots\le |\ln \sqrt{F^*(\varrho)}|\int^\delta_{\sqrt{F^*(\varrho)}} \frac{dt}{\varrho+F(t^2)}\lesssim \frac{\sqrt{F^*(\varrho)}|\ln \sqrt{F^*(\varrho)}|}{\varrho}$$
where the last inequality follows by \eqref{3.2a}. This is the proof of Lemma \ref{l3}.

$\hfill\Box$

Similarly to the proof of \eqref{3.8} we obtain $\displaystyle\int_0^d\frac{\sqrt{F^*(t)}\big|\ln \sqrt{F^*(t)\big|}}{t}dt<\infty$ under hypothesis $\displaystyle\int_0^\delta \big|\ln t \cdot \ln F(t^2)\big|dt<\infty$ for $d,\delta>0$ enough small.\\

Using the general Hardy-Littewood Lemma, we obtain the proof of Theorem \ref{t2}.(ii).

$\hfill\Box$

\section{General Hardy-Littewood Lemma for $f$-H\"older estimates}

We conclude by proving a general Hardy-Littlewood Lemma for $f$-H\"older estimates.
\begin{theorem} \Label{HL}Let $\Om$ be a bounded Lipschitz domain in $\R^N$ and let $\delta_{b\Om}(x)$ denote the distance function from $x$ to the boundary of $\Om$. Let $G:\R^+\to\R^+$ be a increasing function such that $\dfrac{G(t)}{t}$ is decreasing and $\displaystyle\int_0^d\frac{G(t)}{t}dt<\infty$ for $d>0$ small enough.  If $u\in C^1(\Om)$ such that
\begin{eqnarray}\Label{5.1}
\quad |\nabla u(x)|\lesssim \frac{G(\delta_{b\Om}(x))}{\delta_{b\Om}(x)}~~\T{~~for every ~~} x\in \Om.
\end{eqnarray}
Then $|u(x)-u(y)|\lesssim f(|x-y|^{-1})^{-1},\T{~for }x,y\in\Om, x\not=y$ where $f(d^{-1})=\Big(\displaystyle\int_0^d\frac{G(t)}{t}dt\Big)^{-1}.$
\end{theorem}
\begin{remark} If $G(t)=t^\alpha$, Theorem \ref{HL} is the usual Hardy-Littlewood Lemma for domains of finite type. The proof of this theorem in this case can be found in \cite{CS01}. 
\end{remark}
{\it Proof.} Since $u\in C^1$ in the interior of $\Om$, we only need to prove the assertion when $z$ and $w$ are near the boundary. Using a partion of unity, we can assume that $u$ is supported in $U\cap \bar\Om$, where $U$ is a neighborhood of a boundary point $x_o\in b\Om$. After linear change of coordinates, we may assume $x_o=0$ and for some $\delta>0$,
$$U\cap \Om=\{x=(x',x_N)| x_N>\phi(x'), |x'|<\delta, |x_N|\le \delta\},$$
where $\phi(0)=0$ and $\phi$ is some Lipschitz function with Lipschitz constant $M$. Let $x=(x,x_N), y=(y',y_N)\in \Om$ and $d=|x-y|$. For $a\ge 0$, we define the line segment $L_a$ by  
$\theta(x', x_N+a)+(1-\theta)(y',y_N+a)$, $0\le \theta\le 1$. 
Using the Lipschitz property of $\phi$, we obtain
\begin{eqnarray}\Label{5.3}
\begin{split} \tilde x_N+Md=&\theta(x_N+Md)+(1-\theta)(y_N+Md)\\
\ge& Md+\theta \phi(x')+(1-\theta)\phi(y')\\
\ge &Md+\theta (\phi(x')-\phi(\tilde x'))+(1-\theta)(\phi(y')-\phi(\tilde x'))+\phi(\tilde x')\\
\ge& \phi(\tilde x').
\end{split}
\end{eqnarray}  
This implies that $L_a$ lies in $\Om$ for any $a\ge Md$.
Since $u\in C^1(\Om)$,  the Mean Value Theorem tells us that there must exist some $(\tilde x',\tilde x_N+2Md)\in L_{2Md}$ such that 
$$|u(x',x_N+2Md)-u(y',y_N+2Md)|\le |\nabla u(\tilde x',\tilde x_N+2Md)|d.$$
The distance function $\delta_{b\Om}(x',x_N)$ is comparable to $x_N-\phi(x')$, i.e., there are positive constants $c, C$ such that
\begin{eqnarray}\Label{5.2}
c(x_N-\phi(x'))\le \delta_{b\Om}(x',x_N)\le C(x_N-\phi(x')) ~~\T{for }~ x\in \Om\cap U.\end{eqnarray}

Using hypothesis of $G$, combining with \eqref{5.1} and \eqref{5.2}, it follows that 
\begin{eqnarray}
\begin{split}
|u(x',x_N+2Md)-u(y',y_N+2Md)|\lesssim& \frac{G(\delta_{b\Om}(\tilde x',\tilde x_N+2Md))}{\delta_{b\Om}(\tilde x',\tilde x_N+2Md)}d\\
 \lesssim&\frac{G(c(\tilde x_N+2Md-\phi(\tilde x')))}{c(\tilde x_N+2Md-\phi(\tilde x'))}d\\
\lesssim&\frac{G(cMd)}{cMd}d\\ \lesssim &G(d),
\end{split}
\end{eqnarray}
where the last inequality follows by considering two case of $cM$; if $cM<1$, we use $G(t)$ increasing; otherwise, we use $\dfrac{G(t)}{t}$ decreasing. %We consider two case of $cM$. If $cM>1$, then $$
We also have 
\begin{eqnarray}
\begin{split}
|u(x)-u(x',x_N+2Md)|=&\left|\int^{d}_0\frac{\di u(x',x_N+2Mt)}{\di t}dt \right|\\
\lesssim & \int^{d}_0\frac{G(\delta_{b\Om}(x',x_N+2Mt))}{\delta_{b\Om}(x',x_N+2Mt)}dt\\
\lesssim&  \int^{d}_0\frac{G(t)}{t}dt.
\end{split}
\end{eqnarray}

Thus for any $x,y\in \Om$, 
\begin{eqnarray}
\begin{split}
|u(x)-u(y)|\le &|u(x)-u(x',x_N+2Md)|+|u(y)-u(y',y_N+2Md)|\\
&+|u(x',x_N+2Md)-u(y',y_N+2Md)|\\
\lesssim&G(d)+\int^d_0\frac{G(t)}{t}dt\lesssim \int^d_0\frac{G(t)}{t}dt.
\end{split}
\end{eqnarray}
Here, the last inequality follows from
$$G(d)=\int^d_0\frac{G(d)}{d}dt\le\int^d_0\frac{G(t)}{t}dt .$$
This proves the theorem.

$\hfill\Box$

\section*{Appendix}
In this part, we give an explanation of Remark \ref{r1}. First we show the following theorem.
\begin{theorem}\Label{t5} Let $\Om$ and $F$ be defined by \eqref{Om1} or \eqref{Om2} and let $f(d^{-1})=(\sqrt{F^*(d)})^{-1}$ (for $d>0$ small enough). Then $f$-estimate holds for the $\dib$-Neumann problem, that is, 
\begin{eqnarray}
\no{f(\La)u}^2\lesssim \no{\dib u}^2+\no{\dib^* u}^2, 
\end{eqnarray}
 for any $u\in C^\infty_{(0,1)}(\bar\Om)\cap \T{Dom}(\dib^*)$, where $\no{\cdot}$ is the $L^2(\Om)$ norm, $f(\La)$ is the tangential pseudo-differential operator with symbol $f((1+|\xi|^2)^{1/2})$ and $\dib^*$ is the $L^2$-adjoint of $\dib$ with its domain $\T{Dom}(\dib^*)$.
\end{theorem}
{\it Proof.} We will only give the proof in the case $\Om$ is defined by \eqref{Om2}, that is, 
$$\Om\cap U=\{\rho(z)=F(x_1^2)+r(z)<0\},$$
 since the other is proves similarly. Here, $x_1=\Re z_1$.  It is sufficient to show that there exists a family of absolutely bounded weights $\{\Phi^\delta\}$ defined on $S_\delta\cap U$ satisfying
\begin{eqnarray}\Label{5.8}
\underset{i,j=1}{\overset{2}{\sum}}\dfrac{\di^2\Phi^\delta}{\di z_i\di\bar z_j}u_i\bar u_j\simge f(\delta^{-1})^2|u|^2 \quad \T{on }~S_\delta\cap U
\end{eqnarray}
for any $u\in C^\infty_{(0,1)}(\bar\Om\cap U)$, where $S_\delta=\{z\in \Om: -\delta\le \rho(z)\le 0\}$ and $U$ is a neighborhood of the origin (see Theorem 1.4 in \cite{KZ10}).\\

For any $\delta>0$, we define 
$$\Phi^\delta(z):=\exp\left(\dfrac{\rho(z)}{\delta}+1\right)-\exp\left(-\dfrac{x_1^2}{4F^*(\delta)}\right).$$
 The weights $\Phi^\delta$ are absolutely bounded on $S_\delta\cap U$. Computing of the Levi form of $\Phi^\delta$ shows that
\begin{eqnarray}
\begin{split}
\underset{i,j=1}{\overset{2}{\sum}}\dfrac{\di^2\Phi^\delta(z)}{\di z_i\di\bar z_j}u_i\bar u_j
=&\frac{1}{\delta}\left(\underset{i,j=1}{\overset{2}{\sum}}\dfrac{\di^2 \rho(z)}{\di z_i\di\bar z_j}u_i\bar u_j+\frac{1}{\delta}\left|\sum_{j=1}^2 \frac{\di \rho(z)}{\di z_j} u_j\right|^2\right)\exp\left(\dfrac{\rho(z)}{\delta}+1\right)\\
 &+\frac{1}{8F^*(\delta)}\left(1-\frac{x_1^2}{2F^*(\delta)}\right)\exp\left(-\dfrac{x_1^2}{4F^*(\delta)}\right)|u_1|^2\\
\ge&\dfrac{1}{2}\left[\frac{1}{\delta}\frac{ F(x_1^2)}{x_1^2}+\frac{1}{4F^*(\delta)}\left(1-\frac{x_1^2}{2F^*(\delta)}\right)\exp\left(-\dfrac{x_1^2}{4F^*(\delta)}\right)-\frac{c}{F^*(\delta)}\right]|u_1|^2\\
&+cF^*(\delta)^{-1}|u_2|^2,
\end{split}
\end{eqnarray}
for any $z\in S_\delta\cap U$, where $c>0$ will be chosen small.  Here, the inequality follows by the hypothesis of $\rho$ and $F$.\\

We use the notation that
\begin{eqnarray*}
A=\frac{1}{\delta}\frac{ F(x_1^2)}{x_1^2}, \quad B=\frac{1}{4F^*(\delta)}\left(1-\frac{x_1^2}{2F^*(\delta)}\right)\exp\left(-\dfrac{x_1^2}{4F^*(\delta)}\right), \quad\T{and } C=-\frac{c}{F^*(\delta)}.
\end{eqnarray*}
We consider two cases.\\
{\it Case 1.} $x_1^2\le F^*(\delta)$. We have $B\ge \dfrac{e^{-1/4}}{8F^*(\delta)}$, and hence $A+B+C\simge (F^*(\delta))^{-1}$ for a small choice of $c$ in term $C$.\\
{\it Cases 2.} Otherwise, assume $x_1^2\ge F^*(\delta)$. Using our assumption that $\dfrac{F(t)}{t}$ is increasing, we get
$$A=\frac{1}{\delta}\frac{F(x_1^2)}{x_1^2}\ge \frac{1}{\delta}\frac{F(F^*(\delta))}{F^*(\delta)}=\frac{1}{\delta}\frac{\delta}{F^*(\delta)}=\frac{1}{F^*(\delta)}$$
In this case, $B$ can be negative; however, by using the fact that $\min_{t\ge 1/2}\left\{(1-t)e^{-t/2}\right\}=-2e^{-3/2}$ for $t=\dfrac{x_1^2}{2F^*(\delta)}\ge \dfrac{1}{2}$, we have
$B\ge -\dfrac{e^{-3/2}}{2F^*(\delta)}$. This implies $A+B+C\simge (F^*(\delta))^{-1}$ for $c$ small enough\\

Therefore, we obtain \eqref{5.8}. That concludes the proof of Theorem \ref{t5}.

$\hfill\Box$

 By the equivalence of an $f$-estimate on a domain and its boundary in $\C^2$ (see \cite{Kh10}), we have
\begin{eqnarray}\Label{5.13}
\no{f(\La)u}^2\lesssim \no{\dib_bu}^2+\no{\dib_b^* u}^2
\end{eqnarray}
holds for any $u\in C^\infty_{(0,0)}(b\Om)$ or $u\in C^\infty_{(0,1)}(b\Om)$. Here, the norm in  \eqref{5.13} is $L^2$-norm on $b\Om$ and $\dib_b$ is tangential Cauchy-Riemann operator on $b\Om$ with its adjoint $\dib_b^*$.\\

Next, we notice that the hypothesis in Theorem \ref{t2} implies $\underset{t\to 0}{\lim}( t\ln F(t^2))=0$. This limit is equivalent to 
$\underset{\delta\to 0}{\lim}\dfrac{f(\delta^{-1})}{\log \delta}=\infty$. This means superlogarihmic estimates (in the sense of Kohn \cite{Ko02}) for $\Box$ and $\Box_b$ hold. Until now, we did not know if there is a function $f$ such that $f$-H\"older estimate for the integral solution of $\dib$-equation on $\Om$ (defined by \eqref{Om1} or \eqref{Om2}) holds when $F(t^2)=\exp\left(-\dfrac{1}{t|\ln t|^\alpha} \right)$ with $0<\alpha\le 1$. However, in this case, $L^2$-superlogarithmic estimates for $\Box$ and $\Box_b$ still hold.   
\bibliographystyle{alphanum}

\end{document}